%
%
%
%
\input amstex.tex
\input amsppt.sty
\documentstyle{amsppt}
\advance\vsize2\baselineskip
\def\C{\Cal C}
\def\A{\Cal A}

\def\ro{\rho}
\def\bigint{\int}

\def\cal{\Cal}
\def\mathcal{\Cal}
\def\mathbb{\Bbb}
\def\B{\Cal B}

\def\U{\Cal U}

\def\R{\Bbb R}
\def\N{\Bbb N}
\def\T{\Bbb T}

\def\phi{\varphi}
\def\e{\varepsilon}

\def\P{\Cal P}

\def\Id{\operatorname{Id}}

\def\Comp{\cal Comp}

\def\supp{\operatorname{supp}}

\def\min{\operatorname{min}}
\def\max{\operatorname{max}}

\def\Cl{\operatorname{Cl}}

\NoBlackBoxes \topmatter
\title There is no monad based on Hartman-Mycielski functor
\endtitle
\author {Lesya Karchevska$^{1)}$}, {Iryna Peregnyak$^{2)}$} and {Taras Radul$^{3)}$}
\endauthor
\address
$^{1)}$ Department of Mechanics and Mathematics, Lviv National University,
Universytetska st.,1, 79000 Lviv, Ukraine.
\newline e-mail:
 crazymaths\@ukr.net
\endaddress

\address $^{2)}$ Department of Mechanics and Mathematics, Lviv National University,
Universytetska st.,1, 79000 Lviv, Ukraine.
\newline e-mail:irinaperegnyak\@gmail.com
 \endaddress

\address $^{3)}$ Institute of Mathematics, Casimirus the Great University, Bydgoszcz, Poland;\newline
Department of Mechanics and Mathematics, Lviv National University,
Universytetska st.,1, 79000 Lviv, Ukraine.
\newline e-mail:
tarasradul\@yahoo.co.uk
 \endaddress
\leftheadtext {Lesya Karchevska, Iryna Peregnyak and Taras Radul}

\keywords Hartman-Mycielski construction,
 normal functor, monad
\endkeywords
\subjclass 54B30, 57N20
\endsubjclass
\abstract We show that there is no monad based on the normal
functor $H$ introduced earlier by Radul which is a certain
functorial compactification of the Hartman-Mycielski construction
$HM$.
\endabstract
\endtopmatter

\document

\centerline{\bf 0. Introduction}
\vskip 0.3cm

The general theory of functors acting on the category $\Comp$ of
compact Hausdorff spaces (compacta) and continuous mappings was
founded by Shchepin [Sh]. He described some elementary
properties of such functors and defined the notion of the normal
functor which has become very fruitful. The classes  of all normal
 functors include many classical constructions:
the hyperspace {\sl exp}, the space of probability measures $P$, the
space of idempotent measures $I$, and many other functors (cf. [FZ], [TZ], [Z]).

Let $X$ be a space and $d$  an  admissible metric on $X$ bounded
by $1$. By $HM(X)$ we shall denote the space of all  maps from
$[0,1)$ to the space~$X$ such that $f|[t_i,t_{i+1})\equiv{}${\sl const},
for some $0=t_0\le\dots\le t_n=1$, with respect to the  following metric $$
d_{HM}(f,g)=\bigint_0^1d(f(t),g(t))dt,\qquad f,g\in HM(X). $$

The construction of $HM(X)$ is known as the {\sl Hartman-Mycielski
construction} [HM] and was introduced for purposes  of topological groups theory. However it found some applications not connected with groups (see for example [Z1]).

 The construction $HM$ was considered for any compactum $Z$ in [TZ; 2.5.2].
Let $\U$ be the unique uniformity of $Z$. For every $U\in\U$ and
$\e>0$, let
$$<\alpha,U,\e>=\{\beta\in HM_n(Z)\mid m\{t\in
[0,1)\mid (\alpha(t),\beta(t'))\notin U\}<\e\}.$$
The sets
$<\alpha,U,\e>$ form a base of a  topology in
$HM Z$. The construction $HM$ acts also on maps.
Given a map $f:X\to Y$ in $\Comp$, define a map $HM X\to HM Y$
by the formula $HM F(\alpha)=f\circ\alpha$. In general, $HMX$ is not compact.

Let us fix some $n\in\N$. For every compactum $Z$ consider $$\align HM_n(Z)=\Bigl\{f\in
HM(Z)\mid &\text{ there exist } 0=t_1<\dots<t_{n+1}=1\\ &\text{
with }f|[t_i,t_{i+1})\equiv z_i\in Z, i=1,\dots,n\Bigr\}.
\endalign$$
The constructions $HM_n$ define
normal functors in $\Comp$  [TZ; 2.5.2].

 Zarichnyi  has asked if there exists a normal functor in $\Comp$ which contains
 all  functors $HM_n$ as subfunctors (see [TZ]). Such a functor $H$ was constructed in [Ra]. Topological properties of the functor $H$ were investigated in [RR] and [RR1].

 The algebraic aspect of the theory of functors in
categories of topological spaces and continuous maps was inswestigated rather recently. It is based, mainly, on the existence of monad (or triple) structure in the
sense of S.Eilenberg and J.Moore [EM].  We recall the definition of monad  only for the category $\Comp$.
A {\it monad} $\T=(T,\eta,\mu)$ in the category
$\Comp$ consists of an endofunctor $T:{\Comp}\to{\Comp}$ and
natural transformations $\eta:\Id_{\Comp}\to T$ (unity),
$\mu:T^2\to T$ (multiplication) satisfying the relations $\mu\circ
T\eta=\mu\circ\eta T=${\bf 1}$_T$ and $\mu\circ\mu T=\mu\circ
T\mu$. (By $\Id_{\Comp}$ we denote the identity functor on the
category ${\Comp}$ and $T^2$ is the superposition $T\circ T$ of
$T$.)

Many known functors lead to monads: hyperspaces, spaces of
probability measures, superextensions etc. There were many
investigations of monads in  categories of topological spaces and continuous maps(see
for example [RZ] or [TZ]). The following question arises naturally: if the functor $H$ could be completed to a monad?
We give a negative answer in this paper.

\vskip 0.5cm

\centerline{\bf 1. Construction of $H$}

\vskip 0.3cm

Let $X$ be a compactum. By $CX$ we denote the Banach space of all
continuous functions $\phi:X\to\R$ with the usual $\sup$-norm:
$\|\phi\| =\sup\{|\phi(x)|\mid x\in X\}$. We denote the segment
$[0,1]$ by $I$.

For a compactum $X$ let us define the uniformity of $HMX$. For each
$\phi\in C(X)$ and $a,b\in [0,1]$ with $a<b$ we define a
function $\phi_{(a,b)}:HMX\to \R$ by the following formula
$$\phi_{(a,b)}=\frac 1{(b-a)}\bigint_a^b\phi\circ\alpha(t)dt \text{ for some }\alpha\in HMX.$$
Define
$$S_{HM}(X)=\{\phi_{(a,b)}\mid \phi\in C(X) \ \hbox{and} \ (a,b)\subset (0,1)\}.$$

For $\phi_1,\dots,\phi_n\in S_{HM}(X)$ define a pseudometric
$\ro_{\phi_1,\dots,\phi_n}$ on $HMX$ by the formula
$$\ro_{\phi_1,\dots,\phi_n}(f,g)=\max\{|\phi_i(f)-\phi_i(g)|\mid i\in\{1,\dots,n\}\},$$
where $f,g\in HMX$.
The family of pseudometrics
$$\P=\{\ro_{\phi_1,\dots,\phi_n}\mid n\in\N, \ \hbox{where} \
\phi_1,\dots,\phi_n\in S_{HM}(X)\},$$ defines a totally bounded
uniformity $\U_{HMX}$ of $HMX$ (see [Ra]).

For each compactum $X$ we consider the uniform space
$(HX,\U_{HX})$ which is the completion of $(HMX,\U_{HMX})$ and the
topological space $HX$ with the topology induced by the uniformity
$\U_{HX}$. Since $\U_{HMX}$ is totally bounded, the space $HX$ is
compact.

Let $f:X\to Y$ be a continuous map. Define a map $HMf:HMX\to
HMY$ by the formula $HMf(\alpha)=f\circ\alpha$, for all $\alpha\in
HMX$. It was shown in [Ra] that the map
$HMf:(HMX,\U_{HMX})\to(HMY,\U_{HMY})$ is uniformly continuous.
Hence there exists a continuous map $Hf:HX\to HY$ such that
$Hf|HMX=HMf$. It is easy to see that $H:\Comp\to\Comp$ is a
covariant functor and $HM_n$ is a subfunctor of $H$ for each
$n\in\N$.

Let us remark that the family of functions $S_{HM}(X)$ embed $HMX$
in the product of closed intervals $\prod_{\phi_{(a,b)}\in
S_{HM}(X)}I_{\phi_{(a,b)}}$ where $I_{\phi_{(a,b)}}=[\min_{x\in X}
|\phi(x)|,\max_{x\in X}$ $|\phi(x)|]$. Thus, the space $HX$ is the
closure of the image of $HMX$.  We denote by
$p_{\phi_{(a,b)}}:HX\to I_{\phi_{(a,b)}}$ the restriction of the
natural projection.  Let us remark that the function $Hf$ could be
defined by the condition $p_{\phi_{(a,b)}}\circ Hf=p_{(\phi\circ
f)_{(a,b)}}$ for each $\phi_{(a,b)}\in S_{HM}(Y)$.

We will use some properties of the functor $H$ proved in [Ra].
Since the functor $H$ preserves embeddings, we can identify the
space $HA$ with $Hi(HA)\subset HX$ for each closed subset
$A\subset X$ where $i:A\to X$ is the natural embedding. We can
define for each $\alpha\in HX$ the closed set
$\supp\alpha=\cap\{A$ is a closed subset of $X$ such that
$\alpha\in HA\}$. Since $H$ preserves preimages, we have
$\alpha\in H(\supp\alpha)$.

It follows from results of [Ra] and Proposition 5.6 from [Fe] that
there exists a unique natural transformation $\eta:\Id_{\Comp}\to
H$ defined as follows $\eta X(x)(t)=x$ for each $t\in [0,1)$,
where $x\in X$. In other words we have $p_{\phi_{(a,b)}}\circ\eta
X=\varphi$ for each $\varphi\in C(X)$ and $(a,b)\subset(0,1)$.
\vskip 0.5cm

\centerline {\bf 2. Some technical results}

\vskip 0.3cm

It is easy to check that for each $\varphi_1$, $\varphi_2\in
C(X)$, $(a,b)\subset (0,1)$, $\gamma\in HMX\subset HX$ and
$\lambda_1$, $\lambda_2\in \R$ we have
$p_{(\lambda_1\varphi_1+\lambda_2\varphi_2)_{(a,b)}}(\gamma)=\lambda_1p_{\varphi_{1(a,b)}}(\gamma)+\lambda_2p_{\varphi_{2(a,b)}}(\gamma)$.
As well, if $\varphi_1\le\varphi_2$, we obtain
$p_{\varphi_{1(a,b)}}(\gamma)\le p_{\varphi_{2(a,b)}}(\gamma)$.
Since $HMX$ is dense in $HX$, we obtain the following two lemmas.

\proclaim {Lemma 2.1} For each $\varphi_1$, $\varphi_2\in C(X)$, $(a,b)\subset (0,1)$, $\gamma\in HX$ and $\lambda_1$, $\lambda_2\in \R$ we have  $p_{(\lambda_1\varphi_1+\lambda_2\varphi_2)_{(a,b)}}(\gamma)=\lambda_1p_{\varphi_{1(a,b)}}(\gamma)+\lambda_2p_{\varphi_{2(a,b)}}(\gamma)$.
\endproclaim

\proclaim {Lemma 2.2} For each $\varphi_1$, $\varphi_2\in C(X)$, $(a,b)\subset (0,1)$, $\gamma\in HX$ such that $\varphi_1\le\varphi_2$ we have  $p_{\varphi_{1(a,b)}}(\gamma)\le p_{\varphi_{2(a,b)}}(\gamma)$.
\endproclaim

\proclaim {Lemma 2.3} Consider any $\nu\in HX$ and a closed subset $B\subset X$. Then  $\nu\in HB$ iff $p_{\varphi_{1(a,b)}}(\nu)=p_{\varphi_{2(a,b)}}(\nu)$ for each $(a,b)\subset (0,1)$ and
$\varphi_1$, $\varphi_2\in C(X)$ such that $\varphi_1|_B=\varphi_2|_B$.
\endproclaim

\demo {Proof} Necessity. The inclusion $\nu\in HB\subset HX$ means
that there exists $\nu_0\in HB$ with $H(i)(\nu_0) = \nu$, where
$i:B\to X$ is the natural embedding. Hence, for each $(a,b)\subset (0,1)$ and
$\varphi_1$, $\varphi_2\in C(X)$ such that $\varphi_1|_B=\varphi_2|_B$ we have
$p_{\varphi_{1(a,b)}}(\nu)=p_{\varphi_{1}\circ i_{(a,b)}}(\nu_0)=p_{\varphi_{2}\circ i_{(a,b)}}(\nu_0)=p_{\varphi_{2(a,b)}}(\nu)$.

Sufficiency. We can find an embedding $j:B\hookrightarrow
Y$, where $Y\in AR$. Define $Z$ to be the quotient space of the
disjoint union $X\cup Y$ obtained by attaching $X$ and $Y$ by $B$.
Denote by $r:Z\to Y$ a retraction mapping.

Now take any $\nu\in HX\subset HZ$ with the property
$p_{\varphi_{1(a,b)}}(\nu)=p_{\varphi_{2(a,b)}}(\nu)$ for each $(a,b)\subset (0,1)$ and $\varphi_1$,
$\varphi_2\in C(X)$ such that $\varphi_1|_B=\varphi_2|_B$. We
claim that $H(r)(\nu) = \nu$. Indeed, take any $\varphi\in C(Z)$.
Then $p_{\varphi_{(a,b)}}(H(r)(\nu)) = p_{\varphi\circ r_{(a,b)}}(\nu) = p_{\varphi_{(a,b)}}(\nu)$
since $\varphi\circ r|_Y = \varphi|_Y$. Hence, $\nu\in HX\cap HY =
HB$.
\enddemo

Lemmas 2.1 and 2.3 imply the next lemma.

\proclaim {Lemma 2.4} Consider any $\nu\in HX$ and a closed subset $B\subset X$. Then  $\nu\in FB$ iff $p_{\varphi_{(a,b)}}(\nu)=0$ for each $(a,b)\subset (0,1)$ and
 $\varphi\in C(X)$ such that $\varphi|_B\equiv 0$.
\endproclaim

\proclaim {Lemma 2.5} Consider any $\nu\in HX$ and $x\in X$. Then  $x\in\supp\nu$ iff for each neighborhood $O$ of $x$ there exists $a>0$ such that
$p_{\psi_{(0,1)}}(\nu)\ge a$ for each $\psi\in C(X,[0,1])$ such that $\psi|_O\equiv 1$.
\endproclaim

\demo {Proof} Necessity. Suppose the contrary: there exists a neighborhood $O$ of $x$ such that for each   $a>0$ there exists $\psi\in C(X,[0,1])$ such that $\psi|_O\equiv 1$ and $p_{\psi_{(0,1)}}(\nu)< a$. Choose any function $\varphi\in C(X)$ such that $\varphi|_{X\setminus O}\equiv 0$. Let $|\varphi|\le M>0$. By our supposition for each $\varepsilon>0$ we can choose $\psi\in C(X,[0,1])$ such that $\psi|_O\equiv 1$ and $p_{\psi_{(0,1)}}(\nu)< \frac\varepsilon M$. Since $|\varphi|\le M\psi$, we obtain $p_{|\varphi|_{(0,1)}}(\nu)\le p_{M\psi_{(0,1)}}(\nu)<\varepsilon$ using Lemmas 2.1 and 2.2. Thus we have $p_{|\varphi|_{(0,1)}}(\nu)<\varepsilon$ for each $\varepsilon>0$, hence  $p_{|\varphi|_{(0,1)}}(\nu)=0$. It is easy to check that $p_{\varphi_{(0,1)}}(\nu)=0$ too. Then we have $\nu\in H(X\setminus O)$, hence $x\notin\supp\nu$.

Sufficiency. Suppose  $x\notin\supp\nu$. Choose a neighborhood $O$ of $x$ such that $\Cl O\cap\supp\nu=\emptyset$. There exists a function $\psi\in C(X,[0,1])$ such that $\psi|_O\equiv 1$ and $\psi|_{\supp\nu}\equiv 0$. Then $p_{\psi_{(0,1)}}(\nu)=0$ by Lemma 2.4.
 Thus, we obtain a contradiction and the
lemma is proved.
\enddemo

\vskip 0.3cm

\centerline {\bf 3. The main result}

\vskip 0.3cm

For any natural number $n\in\N$ by $K_n$ we denote the finite compactum $\{1,\dots,n\}$ (with discrete topology). Define $\alpha\in HM(K_n\times K_n)\subset H(K_n\times K_n)$ and $\beta\in HM(K_n)\subset H(K_n)$ as follows $\alpha(s)=(i;i)$ and $\beta(s)=i$ if $s\in [\frac{i-1}n,\frac in)$ for $i\in K_n$, $s\in[0,1)$. By $pr_l:K_n\times K_n\to K_n$ for $l\in\{1,2\}$ we denote the natural projections.

\proclaim {Lemma 3.1} We have $(H(pr_1))^{-1}(\beta)\cap(H(pr_2))^{-1}(\beta)=\{\alpha\}$.
\endproclaim

\demo {Proof} Consider any
$\gamma\in(H(pr_1))^{-1}(\beta)\cap(H(pr_2))^{-1}(\beta)$.
Firstly, let us show that $\supp
\gamma\subset\supp\alpha=\{(i;i)|i\in K_n\}$. Suppose the
contrary. The there exist $i,j\in K_n$ such that $i\neq j$ and
$(i;j)\in\supp \gamma$. Consider a function $\psi:K_n\times
K_n\to[0,1]$ such that $\psi(i;j)=1$ and $\psi(k;l)=0$ for each
$(k;l)\neq (i;j)$. By Lemma 2.5 there exists $a>0$ such that
$p_{\psi_{(0,1)}}(\gamma)\ge a$. For $r\in K_n$ define a function
$\varphi_r:K_n\to\R$ by the formula $\varphi_r(s)=1$ if $r=s$ and
$\varphi_r(s)=0$ if $r\neq s$. For $k\in\{1,2\}$ and $r\in K_n$ we
consider the functions $\varphi_r^k=\varphi_r\circ pr_k:K_n\times
K_n\to\R$. Choose a neighborhood $V$ of $\gamma$ defined as
follows $V=\{\gamma'\in H(K_n\times
K_n)\mid|p_{\psi_{(0,1)}}(\gamma)-p_{\psi_{(0,1)}}(\gamma')|<\frac
a2$ and $|p_{\varphi^k_{r(\frac{r-1}n,\frac
rn)}}(\gamma)-p_{\varphi^k_{r(\frac{r-1}n,\frac
rn)}}(\gamma')|<\frac a{2n}$ for each $k\in \{1,2\}$ and $r\in
K_n\}$.

Consider any $\gamma_1\in HMX\cap V$. Since $|p_{\psi_{(0,1)}}(\gamma)-p_{\psi_{(0,1)}}(\gamma_1)|<\frac a2$, we have $m\{t\in[0,1)\mid\gamma_1(t)=(i;j)\}>\frac a2$. Hence there exists $r\in\{1,\dots,n\}$ such that $m\{t\in(\frac{r-1}n,\frac rn)\mid\gamma_1(t)=(i;j)\}>\frac a{2n}$.

If $r\neq i$ we have $p_{\varphi^1_{r(\frac{r-1}n,\frac rn)}}(\gamma_1)=n\bigint_{\frac{r-1}n}^{\frac rn}\varphi^1_r\circ\gamma_1(t)dt<1-\frac a2$.
But $p_{\varphi^1_{r(\frac{r-1}n,\frac rn)}}(\gamma)=p_{\varphi^1_{r}\circ pr_{1(\frac{r-1}n,\frac rn)}}(\gamma)=p_{\varphi_{r(\frac{r-1}n,\frac rn)}}\circ H(pr_1)(\gamma)=p_{\varphi_{r(\frac{r-1}n,\frac rn)}}(\beta)=1$ and we obtain a contradiction with the definition of $V$.

If $r=i$, then we have $r\neq j$ and we obtain a contradiction using similar arguments for the second projection $pr_2$ and the function $\varphi_r^2$.
Hence we have the inclusion $\supp \gamma\subset\supp\alpha$.

Consider any $\varphi\in C(K_n\times K_n)$ and $(a,b)\subset (0,1)$. Define $\psi\in C(K_n)$ as follows $\psi(i)=\varphi(i;i)$ for $i\in K_n$ and put $\xi=\psi\circ pr_1$. Since $\supp \gamma\subset\{(i;i)|i\in K_n\}$, we have $p_{\varphi_{(a,b)}}(\gamma)=p_{\xi_{(a,b)}}(\gamma)$ by Lemma 2.3. Then
$p_{\varphi_{(a,b)}}(\gamma)=p_{\psi\circ pr_1{(a,b)}}(\gamma)=p_{\psi{(a,b)}}(\beta)=p_{\psi\circ pr_1{(a,b)}}(\alpha)=p_{\varphi_{(a,b)}}(\alpha)$. Hence $\alpha=\gamma$.
\enddemo

\proclaim {Theorem 3.2} There is no natural transformation $\mu:H^2\to H$ such that $\mu\circ
H\eta=\mu\circ\eta H=${\bf 1}$_H$.
\endproclaim

\demo {Proof} Suppose that there exists such natural transformation. Let $n\in\N$. For $i\in K_n$ define $\alpha_i\in HM(K_n\times K_n)\subset H(K_n\times K_n)$  as follows $\alpha_i(s)=(i;j)$  if $s\in [\frac{j-1}n,\frac jn)$ for $j\in K_n$ and $s\in[0,1)$. We also define $\A_n\in HM^2(K_n\times K_n)\subset H^2(K_n\times K_n)$  as follows $\A_n(s)=\alpha_i$  if $s\in [\frac{i-1}n,\frac in)$ for $i\in K_n$ and $s\in[0,1)$. Put $H^2(pr_l)(\A_n)=\C_l$ for $l\in\{1,2\}$. Then we have $\C_1=H(\eta K_n)(\beta)$ and $\C_2=\eta HK_n(\beta)$. Hence $\mu K_n(\C_1)= \mu K_n(\C_2)=\beta$. Since $\mu$ is a natural transformation, we have $\mu K_n\times K_n(\A_n)\in (H(pr_1))^{-1}(\beta)\cap(H(pr_2))^{-1}(\beta)$. Hence we obtain $\mu K_n\times K_n(\A_n)=\alpha$ by previous lemma.

By $D$ we denote the two-point set $\{0,1\}$ with discrete topology. For $i\in\{1,\dots,n\}$ define $\gamma_i\in HMD\subset HD$ as follows $\gamma_i(s)=1$  if $s\in [\frac{i-1}n,\frac in)$ and $\gamma_i(s)=0$ otherwise for $s\in[0,1)$. We also define $\B_n\in HM^2D\subset H^2D$ by conditions $\B_n(s)=\gamma_i$  if $s\in [\frac{i-1}n,\frac in)$ for $i\in\{1,\dots,n\}$ and $s\in[0,1)$. Consider a map $f:K_n\times K_n\to D$ defined as follows $f(i;j)=1$ if $i=j$ and $f(i;j)=0$ otherwise. It is easy to see that $Hf(\A_n)=\B_n$. Since $\mu$ is a natural transformation, we have $\mu D(\B_n)=Hf\circ\mu K_n\times K_n(\A_n)=Hf(\alpha)=\eta D(1)$. But it is easy to see that $\B_n$ converges to $\eta HD(\eta D(0))$ if $n\to\infty$. Hence $\mu D$ is not continuous and we obtain a contradiction.
\enddemo


\enddocument